\documentclass[12pt]{amsart}
\usepackage{amssymb,eucal,graphicx}
\usepackage{ifthen}
\usepackage{epsfig}
\usepackage[all,poly,import,color,ps,dvips]{xy}

\SelectTips{cm}{}

\RequirePackage{amsthm}
\RequirePackage{amssymb}
\usepackage[all,ps]{xy}
\usepackage{latexsym}
\numberwithin{equation}{section}

\newtheorem{theorem}[equation]{Theorem}

\newtheorem{corollary}[equation]{Corollary}

\newtheorem{proposition}[equation]{Proposition}

\newtheorem{definition}[equation]{Definition}
\newtheorem{notation}[equation]{Notation}
\newtheorem{example}[equation]{Example}

\newtheorem{remark}[equation]{Remark}

\def\Sing{\operatorname{Sing}}

\def\im{\operatorname{im}}

\def\length{\operatorname{length}}

\def\Card{\operatorname{Card}}

\def\sub{\subseteq}

\def\+{\oplus}                   
\def\*{\otimes}                  

\def\Pic{\operatorname{Pic}}

\def\Supp{\operatorname{Supp}}

\def\Sing{\operatorname{Sing}}
\def\Span{\operatorname{Span}}

\newcommand\Ii{{\mathcal I}}
\newcommand\Oc{{\mathcal O}}

\newcommand\Aa{\mathcal A}

\newcommand\La{\mathcal L}

\newcommand\g{\mathfrak g}
\newcommand\ZZ{\mathbb{Z}}
\newcommand\CC{\mathbb{C}}

\newcommand\Zz{\mathbb Z}
\newcommand\Pp{\mathbb P}

\makeatletter
 \def\ps@notecras{%
    \def\@oddfoot{\small\@date\hfil\slshape\@pauthor\hfil\upshape\thepage}
    \def\@oddfoot{\basdepage}%
   \def\@oddhead{\tetedepage}
  \let\@mkboth\@gobbletwo
 \let\sectionmark\@gobble
\let\subsectionmark\@gobble}
\makeatother
\begin{document}

\title{\bf Singular curves on a $K3$ surface and linear series
on their normalizations}
\author{Flaminio Flamini, Andreas Leopold Knutsen and Gianluca Pacienza}

\address{\hskip -.43cm Flaminio Flamini, Dipartimento di
Matematica, Universit\`a di Roma "Tor Vergata", Via della Ricerca
Scientifica , 1, 00133 Roma, Italy. \linebreak Fax:
+39-06-72594699. e-mail: {\tt flamini@mat.uniroma2.it}.
}

\address{\hskip -.43cm Andreas Leopold Knutsen, Dipartimento di
Matematica, Universit\`a "Roma Tre", Largo S. Leonardo Murialdo,
1, 00146 Roma, Italy. \linebreak Fax: +39-06-54888080. e-mail:
{\tt knutsen@mat.uniroma3.it}. }

\address{\hskip -.43cm Gianluca Pacienza, Institut de Recherche Math\'ematique Avanc\'ee,
Universit\'e \linebreak L. Pasteur et CNRS  7, rue R. Descartes - 67084 Strasbourg Cedex, France.
\linebreak Fax: +33-3-902-40328. e-mail: {\tt pacienza@math.u-strasbg.fr}.  }

\thanks{{\it Mathematics Subject Classification (2000)}: Primary: 14H10; Secondary: 14H51, 14J28, 14J60}
\keywords{curves, $K3$ surfaces, Brill-Noether theory, Severi varieties}

%

\vskip -30 pt

\begin{abstract}

We study the Brill-Noether theory of the normalizations of
singular, irreducible curves on a $K3$ surface. We introduce a
{\it singular} Brill-Noether number $\rho_{sing}$ and show that if
$\Pic(K3)=\Zz [L]$, there are no ${\mathfrak g}^r_d$'s on the
normalizations of irreducible curves in $|L|$, provided that
$\rho_{sing} <0$. We give examples showing the sharpness of this
result. We then focus on the case of {\em hyperelliptic
normalizations}, and classify linear systems $|L|$ containing
irreducible nodal curves with hyperelliptic normalizations, for
$\rho_{sing}<0$, without any assumption on the Picard group.
\end{abstract}

\maketitle

\section*{Introduction}\label{S:intro}

Smooth curves on $K3$ surfaces, and in particular their Brill-Noether theory, have played
a fundamental role in algebraic geometry in the past decades (see e.g.
\cite{MM}, \cite{L}, \cite{GL}, \cite{Muk},  \cite{CLM}, \cite{Fa2}, \cite{V1},  \cite{AV}, \cite{A},
\cite{FP} and \cite{V1bis}, just to mention a few). The Brill-Noether theory of these curves is both an important subject
in its own right, especially because it is connected to the geometry of the surface and, at the same time, it is
an important tool to prove results about smooth curves with general moduli.

Recall the following fundamental result of R.~Lazarsfeld:

\medskip

\noindent {\bf  Theorem (Lazarsfeld \cite{L}).} {\it Let $S$ be a
$K3$ surface with $\Pic (S)= {\mathbb Z}[L]$, with $L^2=2p-2 >0$.
Let $\rho(p,r,d):=p-(r+1)(p-d+r)$ be the {\em Brill-Noether
number}.

\begin{itemize}
\item[(i)] If $X\in |L|$
 is smooth, and $\rho(p,r,d) < 0$, then
$$
 W^r_d (X):= \{ \Aa \in \; \Pic^d(X) \; | \; h^0(X, \Aa) \geq r+1 \}= \emptyset,
$$
\item[(ii)] If $X\in |L|$ is a general member and $\rho(p,r,d)\geq 0$, then $W^r_d (X)$ is smooth
outside $W^{r+1}_d (X)$, and of
dimension $\rho(p,r,d)$.
\end{itemize}
}

\medskip

Similarly, in \cite{GL}, Green and Lazarsfeld proved that all the
smooth curves in a complete, base point free linear system $|L|$
on a $K3$ surface have the same Clifford index. Moreover, they
proved that if the Clifford index is non-general, i.e. less than $
\lfloor \frac{p-1}{2} \rfloor $, where $p=p_a(L)$ is the
arithmetic genus of the curves in $|L|$, then it is induced by a
line bundle on the surface. As a consequence,  one can easily
construct curves of any prescribed Clifford index and gonality on
$K3$ surfaces \cite{K2}. This shows that curves on $K3$ surfaces
with arbitrary Picard group may be far from being generic with
respect to the Brill-Noether theoretical point of view. The
hyperelliptic case is ``classical'' and dates back to Saint-Donat
\cite{SD}, who classified linear systems on $S$ containing
hyperelliptic curves :

\medskip

\noindent
{\bf Theorem (Saint-Donat \cite{SD}).} {\it Let $L$ be a globally generated line bundle with $L^2 >0$ on a
$K3$ surface $S$. The following conditions are equivalent:
\begin{itemize}
\item[(i)]  the morphism $\varphi_L$  defined by $|L|$, is not birational.
\item[(ii)]  There is a smooth, hyperelliptic curve in $|L|$.
\item[(iii)]  All the smooth curves in $|L|$ are hyperelliptic.
\item[(iv)] Either $L^2 = 2$, or there exists a smooth, irreducible, elliptic curve $E$ on $S$ satisfying
$E.L=2$, or $L \sim 2B$ for a smooth, irreducible curve $B$ with $B^2=2$.
 \end{itemize}
A linear system $|L|$ satisfying these properties is said to be {\em hyperelliptic}.
}

\medskip

Whereas the Brill-Noether theory of {\it smooth} curves on $K3$ surfaces is quite well understood,
almost nothing is known about the Brill-Noether theory of {\it singular} curves; by this we mean
the Brill-Noether theory of {\it normalizations} of singular curves on $K3$ surfaces.
Singular curves on $K3$ surfaces with nodes as the only singularities
behave in the expected way in the following sense, as a consequence of a result of Xi Chen and of the theory
of Severi varieties (cf. Remark \ref{rem:sevvar} for details):

\medskip

\noindent
{\bf Theorem (Xi Chen \cite{C}).} {\it If
$S \subset {\mathbb P}^p$, $ p \geq 3$, is a general, primitive
$K3$ surface of genus $p$, whose hyperplane divisor is $H$,
then, for each integer $m \geq 1$, the Severi varieties
$V_{|mH|, \delta}(S) \subset |mH|$, parametrizing
the universal family of irreducible curves in $|mH|$ having
$\delta$ nodes as the only singularities, are nonempty and hence of the
expected dimension $p_a (mH) -\delta$, for each $\delta \leq p_a(mH)$.}

 \medskip

It is therefore natural to try to extend the study of linear series on smooth curves lying
on a $K3$ surface to the case of curves which are the normalizations of singular curves
on a $K3$ surface, at least in the nodal case.
Beside of the interest on its own, the problem and the results
we obtain in this paper are
related to the study of the Hilbert scheme of
points $S^{[k]}$ of a $K3$ surface. Indeed, a $\g^1_{k}$ on the
normalization $C$ of a curve $X$ on $S$ gives rise to a rational
curve in $S^{[k]}$
and, as a consequence of a result of D. Huybrechts and S. Boucksom
(\cite[Prop. 3.2]{H2} and \cite[Th. 1.2]{Bou}),
we know that rational curves span the Mori cone
of $S^{[k]}$. (When $k=2$, see \cite{HT} for a conjecture about the
numerical and geometrical properties of the  rational curves generating
the extremal rays.) We plan to return to this topic in a future work.

In this paper we draw the attention to a number which seems to be  an important invariant
for the study of the Brill-Noether theory of normalizations of singular curves on $K3$ surfaces.

Let $X$ be a reduced, irreducible, possibly singular curve in
$|L|$ and let $g$ be its {\em geometric genus} (i.e. the genus of
its normalization). Set $p_a := p_a(X)$, and define the {\it
singular} Brill-Noether number as follows :
\begin{equation}\label{eq:(*)}
\rho_{sing}(p_a,r,d,g):= \rho(g,r,d)+ p_a - g.
\end{equation}
Our starting point is a result which we derive from the
interesting work of  T. L. G\'omez \cite{G} concerning rank-one,
torsion-free sheaves on singular curves on smooth $K3$ surfaces.

\medskip

\noindent {\bf Theorem 1}. {\em Let $S$ be a $K3$ surface such
that $\Pic (S) = {\mathbb Z}[L]$ with $L^2 > 0$. Let $X$ be a
reduced, irreducible curve in the linear system $|L|$, of
arithmetic genus $p_a := p_a(X) = p_a(L)$ and having $\delta$
singular points of given multiplicities as the only singularities.
Let $C\stackrel{\nu}{\to} X$ be the normalization of $X$, whose
geometric genus is denoted by $g$. Then, if
$\rho_{sing}(p_a,r,d,g) <0$ one has $W^r_d (C)= \emptyset$. }

\medskip

Note that, when $X$ has $\delta$ nodes, then $p_a - g = \delta$.
In particular, for $r=1$ and $d=2$, we have $$
 \rho_{sing}(p_a,1,2, p_a - \delta) < 0 \Leftrightarrow \delta \leq
(p_a(X) - 3)/2.
$$

We show the sharpness of Theorem 1 by
providing examples (cf. Examples \ref{exm1} and \ref{exm2}) of
nodal curves on a $K3$ surface having hyperelliptic normalizations, with
$\rho_{sing}\geq 0$.

Examples \ref{exm1} and \ref{exm2} also show that the theorem of
Green and Lazarsfeld on the constancy of the Clifford index of
{\it smooth} curves in a linear system $|L|$ cannot be extended to
the curves in the Severi varieties  $V_{|L|, \delta}(S)$ of $|L|$
(cf. Remark \ref{rem:exa}).

Our second main result is a classification, via a Bogomolov-Reider
type approach, of linear systems $|L|$ on a $K3$ surface $S$ with
{\em arbitrary} Picard group, containing $\delta$-nodal curves $X
\in |L|$ with hyperelliptic normalizations and with $\delta \leq
(p_a (X)-3)/2$ (that is $\rho_{sing}(p_a(X),1,2,p_a(X) - \delta)
<0$).

In order to achieve the classification, we have to recall the
classical notion of neutral node for
irreducible, $\delta$-nodal curves having
hyperelliptic normalizations (cf. \S \,\ref{S:nodchar}). More
precisely, if $C$ is the normalization of such a curve $X$ and if
$y_1, y_2 \in C$ are the preimages of a node $x \in X$ then, if
$y_1$, $y_2$ are conjugated by the hyperelliptic involution
existing on $C$, $x \in X$ is called a {\em neutral node},
otherwise $x \in X$ is a {\em non-neutral node} (cf. \,Definition
\ref{def:neutral}).

Our result shows that any element in $V_{|L|, \delta}$ with hyperelliptic
normalization cannot have too few non-neutral nodes:

\medskip

\noindent {\bf Theorem 2}. {\em Let $S$ be a $K3$ surface and $L
\in \Pic (S)$ be globally generated and nonhyperelliptic.
Let $\delta$ be a positive integer such that
$p_a(L) - \delta \geq 2$.

Inside the {\rm Severi variety} $V_{|L|, \delta}$, consider
the subscheme $V_{|L|, \delta}^{hyper} $ parametrizing those
elements having  hyperelliptic normalizations. Let $[X] \in
V_{|L|, \delta}^{hyper}$ and let $\delta_0 \leq \delta$ be the number
of non-neutral nodes of $X$.

Then either $\delta_0 \geq (p_a(L) -5) / 2$ or $(\delta_0, p_a(L))=(2,10)$.

Moreover, if  $\delta_0 \leq (p_a(L) -3) / 2$ or $(\delta_0,
p_a(L)) = (2,10)$ then, for any pair of (possibly coinciding)
points on $X$ belonging to the image of the $\g^1_2$ of $C$, there
exists an effective divisor $D$ on $S$ passing through all the
non-neutral nodes of $X$ and through the pair of points and
satisfying one of the following conditions:

\begin{itemize}
\item[(i)] $L = {\mathcal O}_S(3 D),\ \delta_0 = 2,\ p_a(L) = 10$;
\item[(ii)] $L ={\mathcal O}_S(2 D)$ and either $\delta_0 =
(p_a(L) -3) / 2$ and $D^2 = \delta_0+1 \geq 4$ (so $\delta_0 \geq
3$ is odd), or $\delta_0 = (p_a(L) -5) / 2$ and $D^2 = \delta_0+2$
(so $\delta_0$ is even); \item[(iii)] $D$ and $X$ are linearly
independent in $\Pic (S)$, $D.X=2\delta_0 +2$, $(p_a(L) -5) / 2
\leq \delta_0 \leq (p_a(L) -3) / 2$ and $D^2 =\delta_0$ (so
$\delta_0$ is even).
\end{itemize}
}

\medskip

A similar statement can also be obtained for curves with
normalizations possessing $\g^1_k$'s for $k >2$, but to avoid
longer classification results and proofs - which are more
difficult only from a technical point of view - we choose to keep
the exposition clearer and state it only for $k=2$.

Theorem 2 has the following two corollaries:

\noindent {\bf Corollary 3}. {\em Let $L$ be a globally generated, nonhyperelliptic
line bundle on a $K3$ surface and $\delta$ be a positive integer such that
$p_a(L) - \delta \geq 2$.

If $V_{|L|, \delta}^{hyper} \neq \emptyset$, then either
$\delta \geq (p_a(L) -5) / 2$ or $(\delta, p_a(L))=(2,10)$ }

\medskip

\noindent {\bf Corollary 4}. {\em Let $S$ be a $K3$ surface with $\Pic (S) = {\mathbb Z} [L]$, for
some line bundle $L$ on $S$ with $L^2 \geq 2$, and $n \geq 1$ an
integer.

If $V_{|nL|, \delta}^{hyper}\neq \emptyset$, for some $\delta \leq
\frac{p_a(nL)-3}{2}=\frac{1}{4}n^2L^2-1$ then either
\begin{itemize}
\item[(a)] $(n,L^2)=(3,2)$ and $\delta=2$ or $3$, or
\item[(b)] $\sqrt{\frac{\delta + 1}{L^2}} \in \ZZ$ and $n=2\sqrt{\frac{\delta
+ 1}{L^2}}$, or
\item[(c)] $\sqrt{\frac{\delta + 2}{L^2}} \in \ZZ$
and $n=2\sqrt{\frac{\delta + 2}{L^2}}$.
\end{itemize}
In particular, in any case we have $n \geq 2$.}

\vskip 10pt

The paper is organized as follows. After some preliminaries, we
prove Theorem 1 in Section 1. Section 2 is devoted to discussing
two examples that prove the sharpness of Theorem 1. In Section 3
we briefly recall the notion of neutral and non-neutral nodes and
give some results on the matter. In the last section, we translate the
fact of having hyperelliptic normalization into a non-separation
statement for $|L|$ and then apply a Bogomolov-Reider type
technique to prove Theorem 2 and its two corollaries.

\section*{Preliminaries}\label{S:pre}
%

We work in the category of algebraic
$\CC$-schemes. Terminology and notation are the standard ones
(cf. \cite{Ha}).

By a $K3$ surface is meant a smooth, algebraic surface $S$ (reduced and irreducible)
with $h^1(\Oc_S)=0$ and $K_S=0$. The {\it sectional genus} $p_a(L)$ of a line bundle $L$
(with sections) is defined by $L^2=2p_a(L)-2$. By adjunction, $p_a(L)$ is the arithmetic genus
of any member of $|L|$. Recall that by Saint-Donat's classical results \cite{SD}, on a $K3$ surface $|L|$
is base point free if and only if it has no fixed components.

We use the following standard definitions:

\begin{definition}\label{def:prim}
A line bundle $L$ on a surface $S$ is called {\rm primitive} if $L
\not \sim n L' $ for any $n >1$ and $L' \in \Pic (S)$.

A {\rm polarized surface} (resp. {\rm primitively polarized
surface}) is a pair $(S, L)$, where $S$ is a surface and $L \in
\Pic( S)$ is globally generated and ample (resp. primitive,
globally generated and ample).
\end{definition}

Let $S$ be a smooth, algebraic surface and let
$|L|$ be a complete linear system on $S$, whose general element
is a smooth, irreducible curve. For any integer
$0 \leq \delta \leq p_a(L)$, one denotes
by
\begin{equation}\label{eq:sevvar}
V_{|L|, \delta}(S)
\end{equation}the locally closed and functorially defined
subscheme of $|L|$ parametrizing the universal family of irreducible curves in
$|L|$ having $\delta$ nodes as the only singularities. These
are classically called {\em Severi varieties} of irreducible, $\delta$-nodal
curves in $|L|$ on $S$.

\begin{remark}\label{rem:sevvar}\normalfont{
Recall that if $S$ is a $K3$ surface, there are several properties
concerning Severi varieties on $S$ which are well-known.
Indeed, let $L$ be a globally generated
line bundle on $S$ such that $p_a(L) \geq 2$. Assume that, for a given
$\delta \leq p_a(L)$, $V_{|L|, \delta}(S) \neq \emptyset$.
Let $[X] \in V_{|L|, \delta}(S)$ and let $N = \Sing(X)$ be the subscheme of its nodes.
Since $S$ is a $K3$ surface, then $N$ always
imposes independent conditions on the linear system $|L|$, i.e.
$H^1(S, \Ii_{N/S} \otimes L) = (0)$ (see e.g. \cite{CS}). By the description of the
tangent space at $[X]$ to $V_{|L|, \delta}(S)$, the condition $H^1(S, \Ii_{N/S} \otimes L) = (0)$
is equivalent to the fact that $[X]$ is an unobstructed point of
$V_{|L|, \delta}(S)$ (see e.g. 
\cite{CS} or \cite{F}).

Therefore, on $K3$ surfaces, if $V_{|L|, \delta} (S) \neq \emptyset$, then
each of its irreducible components is {\em regular}, i.e.
smooth and of the expected dimension $p_a(L) - \delta$.
Furthermore, the non-emptyness of $V_{|L|, \delta}(S)$, for some $\delta$,
implies, for each $\delta' < \delta$, that $V_{|L|, \delta'}(S) \neq \emptyset$ and, moreover,
that each irreducible component of $V_{|L|, \delta'}(S)$ is also regular.

For what concerns existence results on $K3$ surfaces we refer to Chen's result
(cf. \cite{C}) and its consequence quoted in the introduction.
On the other hand, we stress that it is still not known whether
the Severi varieties  $V_{|mH|, \delta}(S) $ are irreducible.}
\end{remark}

\begin{definition}\label{def:sevhyp}
Let $S$ be a smooth, algebraic surface and let $|L|$ be a complete
linear system on $S$, whose general element is smooth and irreducible.
Consider $V_{|L|, \delta}(S)$. We denote by $$V_{|L|, \delta}^{hyper}(S)$$the subscheme of
$V_{|L|, \delta}(S)$ parametrizing curves whose normalizations
are 
hyperelliptic curves.
Curves in $V_{|L|, \delta}^{hyper}(S)$ will be called $\delta$-{\em hyperelliptic nodal curves}.
\end{definition}

\noindent
When no confusion arises, we will omit $S$ and simply write
$V_{|L|, \delta}$ and $V_{|L|, \delta}^{hyper}$.

%
\section {An extension of Lazarsfeld's result to the singular case}\label{S:1}
%

The aim of this section is to prove Theorem 1 stated in the introduction.

\begin{remark}\label{rem:Laz}\normalfont{
 Before proving the theorem, we want to add some comments on the
hypothesis $ \Pic(S) = {\mathbb Z}[L]$. This could seem
restrictive, but one cannot hope to weaken this hypothesis too
much if one uses the Lazarsfeld and the G\'omez procedure (cf. the
final Remark in \cite{L}). However, going through the paper of
Lazarsfeld, one sees that the crucial point is to make sure that a
certain vector bundle has no nontrivial endomorphisms (cf.
\cite[Prop. 2.2]{L}). This is the case under the following
hypothesis:
\begin{equation}
  \label{eq:suff}
 L \; \mbox{has no decomposition} \; L \sim M+N \; \mbox{with} \;
h^0(M) \geq 2 \; \mbox{and} \;  h^0(N) \geq 2,
\end{equation}
if one carefully looks through the proof of \cite[Lemma 1.3]{L}.

Therefore, the hypothesis $ \Pic(S) = {\mathbb Z}[L]$ in Theorem
1, as well in Lazarsfeld's theorem quoted in the introduction, can
be weakened to \eqref{eq:suff}.}
\end{remark}

To prove Theorem 1
we first need to introduce some
notation and some more general facts.

Let $Y$ be an integral curve (not necessarily smooth), of
arithmetic genus $p_a$; one can consider the complete scheme
$\overline{J}^d(Y)$, defined as the scheme parametrizing rank one
torsion-free sheaves on $Y$ of degree $d$ (recall that, if $\Aa
\in \overline{J}^d(Y)$, by definition $d= \deg(\Aa) := \chi(\Aa) +
p_a -1$, cf. e.g. \cite[\S\ 2]{G}). When $Y$ lies on a smooth
surface, $\overline{J}^d(Y)$ is integral and contains, as an open,
dense subscheme, the generalized Jacobian $J^d(Y)$, defined as the
variety parametrizing line bundles of degree $d$ on $Y$. Therefore
$\overline{J}^d(Y)$ is a compactification of $J^d(Y)$. One can
consider the complete subscheme of $\overline{J}^d(Y)$, called the
{\em generalized Brill-Noether locus},
\begin{equation}\label{eq:aio1}
\overline{W}^r_d(Y)
\end{equation}parametrizing those
sheaves $\Aa$ which are special and such that $h^0(\Aa) \geq r+1$.

\begin{proposition}\label{prop:ballfo}
Let $S$ be a $K3$ surface and $L$ be a globally generated line
bundle on $S$ such that $L^2 > 0$.  Let $X$ be a reduced,
irreducible curve in $|L|$, of arithmetic genus $p_a := p_a(X) =
p_a(L)$ and having $\delta$ singular points of given
multiplicities as the only singularities.

Let $C\stackrel{\nu}{\to} X$ be the normalization of $X$, whose
geometric genus is denoted by $g$.

Assume furthermore that $C$ possesses a line bundle $\Aa \in
\Pic^d(C)$ such that $|\Aa| = \g^r_d$; then $\nu_*(\Aa)
\in\overline{W}^r_{d + p_a - g}(X).$
\end{proposition}

\begin{proof} It is well-known that $\nu_*(\Aa)$ is a
rank-one, torsion-free, coherent sheaf on $X$ (see e.g. \cite[II, 5.20]{Ha}). Therefore
we have that
\begin{equation}\label{eq:degree1}
\deg(\nu_*(\Aa)) = \chi(X, \nu_*(\Aa)) + p_a(X) -1.
\end{equation}
By Leray's isomorphism,
\begin{equation}\label{eq:degree2}
\chi(X,\nu_*(\Aa)) = \chi(C,\Aa) = d - g +1 \; \; {\rm and} \;\;
h^0(X,\nu_*(\Aa))= h^0(C, \Aa)= r+1. \end{equation} This means
that $\deg(\nu_*(\Aa)) = d-g+1 + p_a(X) - 1 =  d + (p_a(X) - g)$.
\end{proof}

By recalling Formula \eqref{eq:(*)} in the introduction, we can
state the following:

\begin{proposition}\label{cl:aio6}
With notation and assumptions as in Theorem 1, 
we have:
\begin{equation}\label{eq:claio6}
\overline{W}^r_{d+ p_a - g}(X) = \emptyset \Leftrightarrow
\rho_{sing}(p_a,r,d,g) <0.
\end{equation}
\end{proposition}

\begin{proof} ($\Rightarrow$) Assume $\rho_{sing}(p_a,r,d,g) \geq
0$. By the very definition of $\rho_{sing}$ and by a
straightforward computation, we have that
\begin{equation}\label{eq:referee}
\rho_{sing}(p_a,r,d,g) = \rho (p_a, r, d + p_a - g).
\end{equation} Since $p_a = p_a(X) = p_a(L)$, thus for the general $D\in |L|$ we have
$\dim(W^r_{d+m}(D)) = \rho(p_a,r,d + p_a - g)$ by Lazarsfeld's
theorem quoted in the introduction. By upper semi-continuity, as
in \cite[Proof of Theorem I, p.~749]{G}, we have that
$\overline{W}^r_{d+ p_a - g}(X)$ is nonempty and
$\dim(\overline{W}^r_{d+ p_a - g}(X) ) \geq \rho(p_a,r,d + p_a -
g) \geq 0$.

\noindent ($\Leftarrow$) Assume $\overline{W}^r_{d+ p_a - g}(X)
\neq \emptyset$; by \cite[Prop. 2.5 and proof of Thm. I at
p.~749-750]{G} it is possible to deform the couple $(X, \Aa)$,
where $\Aa\in\overline{W}^r_{d+ p_a - g}(X)$, to a couple $(D,
\La)$, where $D$ is  a smooth, irreducible curve in $|L|$ and $\La
\in W^r_{d+ p_a - g}(D)$. By Lazarsfeld's theorem quoted in the
introduction, since $p_a = p_a(X) = p_a(L)$ and $D$ is smooth,
then $\rho(p_a,r,d + p_a - g) \geq 0$. One can conclude by using
the equality in \eqref{eq:referee}.
\end{proof}

\begin{proof}[Proof of Theorem 1.] It directly follows from Propositions \ref{prop:ballfo} and
\ref{cl:aio6}.
\end{proof}

Suppose now that $X$ is a reduced, irreducible curve on $S$ having
$\delta$ nodes as the only singularities. Recall that in this case
we have that $p_{a}(X) - g= \delta$, since each singular point is
an ordinary double point. In particular, when $r=1$ and $d=2$, the
condition  $\rho_{sing}(p_a(X),1,2,p_a(X) - \delta)<0$ is
equivalent to the upper-bound on the number of nodes
$
 \delta \leq (p_a (X)-3)/2.
$

 Thus, recalling Definition \ref{def:sevhyp}, an immediate
 consequence of  Theorem 1
 is the following.

\begin{corollary}\label{cor:vhdh}
Let $S$ be a $K3$ surface such that $\Pic(S) = \ZZ[L]$, with $L$
of sectional genus $p=p_a(L) \geq 3$. If
\begin{equation}\label{eq:bound}
\delta \leq (p-3)/2
\end{equation}then $V_{|L|, \delta}^{hyper} = \emptyset$.
\end{corollary}

\begin{remark}\label{rem:notfine}\normalfont{To conclude this section,
we remark that even if non-existence results easily follow from
Lazarsfeld's and G\'omez's procedures, the same does not occur for
what concerns existence results in the complementary range of
values $\rho_{sing}(p_a,r,d,p_a - \delta) \geq 0$. To see this,
suppose for simplicity that $X$ is nodal. Let $N$ be the scheme of
its $\delta$ nodes. Let $W^r_{d}(X,N)\subset \overline{W}^r_d(X)$
be the subscheme parametrizing the elements of
$\overline{W}^r_d(X)$ which are not locally free exactly at $N$.
Then $W^r_{d}(X,N)\cong W^r_{d-\delta}(C)$ (cf. e.g. \cite{BF}).

Proposition \ref{cl:aio6} yields that, if $\rho_{sing}\geq 0$,
then $\overline{W}^r_d(X)\not=\emptyset$, but this does not insure
that $W^r_{d}(X,N)\subset \overline{W}^r_d(X)$ is nonempty. We
shall discuss in \S\,\ref{S:2} two examples of existence when
$\rho_{sing}\geq 0$. These examples will show the sharpness of
Theorem 1. }
\end{remark}

%
\section {Existence of $\delta$-hyperelliptic nodal curves with $\rho_{sing}\geq 0$}\label{S:2}
%
In this section $X$ will always denote a reduced, irreducible curve
on $S$ having $\delta$ nodes as the only singularities.

Let $S$ be as in Corollary \ref{cor:vhdh} and suppose further that
it is a general, primitive $K3$. Then, by
Chen's theorem quoted in the introduction,
$V_{|L|, \delta} \neq \emptyset$, for each $\delta \leq p$. On the
other hand, \eqref{eq:bound} gives an upper-bound on $\delta$,
which ensures that the subschemes $V_{|L|, \delta}^{hyper}$ are
empty.
However, these subschemes are not always empty. Indeed, for
$\delta = p-2$, one clearly has
\[
V_{|L|, p-2} = V_{|L|, p-2}^{hyper},
\]
since the normalizations of the curves
they parametrize have genus $2$. These can be considered as trivial examples of
$\delta$-hyperelliptic nodal curves (cf.\ Definition \ref{def:sevhyp}).
Therefore nontrivial questions on $V_{|L|, \delta}^{hyper}$ make
sense only for
\begin{equation}\label{eq:nontriv}
p \geq 3 \; \; {\rm and} \; \; \frac{p-2}{2} \leq \delta \leq p -3.
\end{equation}

We shall now discuss in details the unique "nontrivial" examples of
$\delta$-hyperelliptic nodal curves on a
general complete intersection $K3$ surface of type $(2,3)$ in $\Pp^4$
and of type $(2,2,2)$ in $\Pp^5$ (see Examples \ref{exm1} and
\ref{exm2} below). We start with the following remark and result.

\begin{remark} \label{rem:immersione}
 \normalfont{Recall from \cite{SD} that if $L$ is a globally generated line bundle on a $K3$ surface
$S$ with $p:=p_a(L) \geq 3$ and such that $|L|$ is
non-hyperelliptic, i.e. the morphism $\varphi_L$ defined by $|L|$
is birational (cf. Saint-Donat's Theorem quoted in the
introduction). Thus, $\varphi_L$ maps $S$ birationally onto a
projectively normal surface $\varphi_L(S) \sub \Pp^p$ with at
worst rational double points as singularities. These points are
exactly the images of the contracted, smooth rational curves on
$S$, i.e. curves $\Gamma$ such that $\Gamma.L=0$. Moreover
$\varphi_L$ is an isomorphism outside these contracted curves.

It follows that every irreducible curve $X \in |L|$ is mapped isomorphically
to a hyperplane section of $\varphi_L(S)$.}
\end{remark}

\begin{proposition}\label{prop:projgeom}
Let $S$ be a $K3$ surface and $L$ be a globally generated line
bundle on $S$ with $p:=p_a(L) \geq 3$ and such that $|L|$ is
non-hyperelliptic. Let $[X] \in V_{|L|, \delta}$, for some $\delta
\leq p-3$,  and let $N= \Sing(X)$. Let $\nu:C \to X$ be its
normalization, where $g=p-\delta$ is the geometric genus of $C$.
Identify $X$ with its image in $\varphi_L(S) \subset \Pp^p$ by
Remark \ref{rem:immersione}.

Then $C$ is hyperelliptic if and only if   the rational map
$$
 {\pi_{\Lambda_N}}_{|X} : X - - \to \Pp^{g-1},
$$ induced on $X \subset \Pp^{p-1}$ by the projection  from the
linear subspace $\Lambda_N:=\Span_{\Pp^{p-1}}(N)$, maps $X$ $2:1$
onto a smooth, rational normal curve in $\Pp^{g-1}$.
\end{proposition}

\begin{proof} Let $\Delta_C$ be the degree-two pencil on
 $C$. Then, $\omega_C \cong \Delta_C^{\otimes (g-1)}$,
i.e. the canonical morphism
$$
\Phi_C : C \to \Pp^{g-1}
$$
maps $C$ $2:1$ onto a rational normal curve $\Gamma \subset \Pp^{g-1}$,
where $\Phi_C$ is the morphism given by $|\omega_C|$.

We have $X \cong \varphi_L(X)=\varphi_L(S) \cap H$, for some
$H=\Pp^{p-1}$. Consider $\Lambda_N $ the linear subspace of $H$
given by the span of the $\delta$ nodes of $X \subset H$. By our
assumption ${\rm Card}(N) = \delta \leq p-3$, we have that
$\dim(\Lambda_N) = \delta -1$ if and only if the $\delta$ nodes
are in general linear position in $H \cong \Pp^{p-1}$, i.e. if and
only if $h^1(\Ii_{N/S}(L)) = 0$. This vanishing holds by Remark
\ref{rem:sevvar}.

Therefore, for each $\delta \leq p-3$, we have
\begin{equation}\label{eq:9lu2}
\dim (\Lambda_N) = \delta-1.
\end{equation}

Let $\tilde{N}$ be the $2\delta$ points given by the pull-back of
$N$ to $C$. Since $C$ is the normalization of $X$, by adjunction
on the surface $\tilde{S}$ (given by the blow-up of $S$ along
$N$), we know that
\begin{equation}\label{eq:9lu1}
\omega_C \cong \Oc_C(\nu^*(L)(-\tilde{N}));
\end{equation}
thus,
\begin{equation}\label{eq:9lu*}
 \nu_*(\omega_C) \cong \Ii_{N/X}(L),
\end{equation}
which gives
the following commutative diagram:
\begin{equation} \label{eq:9lu3}
 \xymatrix{
     C  \ar[dr]_{\Phi_C} \ar[r]^{\nu} & X \ar[d]^{\tau_X} & \hspace{-0.8cm} = \varphi_L(S) \cap H \sub \Pp^{p-1}  \\
     & \Gamma & \hspace{-3cm} \sub \Pp^{g-1},
 }
\end{equation}
where $\tau_X$ is the map induced on $X$ by $|\Ii_{N/X}(L)|$, that
is the one induced on $X$ by the linear projection from
$\Lambda_N$ $$
 \pi_{\Lambda_N}: \Pp^{p-1} - - - \to \Pp^{g-1}.
$$
Therefore, by Diagram \eqref{eq:9lu3} we can conclude.
\end{proof}

The previous proposition can be rephrased as follows: $C$ is hyperelliptic
if and only if
the projection from the span of the nodes of $X$ maps $X$ onto a
smooth curve in $\Pp^{g-1}$ lying on a number of hyperquadrics
which is bigger than the expected one.

From \eqref{eq:nontriv}
it follows that, when $(S,H)$ is primitive, and $H$ is very ample (that is, $p_a(H) \geq 3$),
the first nontrivial example to
consider is the following.

\vskip 7pt

\begin{example} \label{exm1}
\normalfont{Let $S = Q_2 \cap Q_3 \subset {\mathbb P}^4$ be a
general, primitive, smooth $K3$ surface which is the complete
intersection of a quadric and a cubic hypersurface. If $|H|$
denotes the hyperplane linear system of $S$, such that $\Pic(S) =
\ZZ[H]$, then $p_a(H) =4$. From  \eqref{eq:nontriv}, we want to
consider $\delta = 1 = p_a(H) - 3$; by Remark \ref{rem:sevvar}, we
know that $V_{|H|,1} \neq \emptyset$ on such a $S$.

Take $q \in S$ a general point and
consider $H_q := T_q (Q_2)$ the hyperplane in ${\mathbb P}^4$ which is
the tangent space to $Q_2$ at the point $q$. Let $Q_q := Q_2 \cap
H_q$, which is a quadric cone in $H_q \cong {\mathbb P}^3$. Thus, the
curve $X: = H_q \cap S$ is such that $[X] \in V_{|H|,1}$ on $S$.
Moreover $\deg (X)=6$ and its normalization  $C$ has genus $g = 3$.

Observe that, since $Q_q$ is a quadric cone in ${\mathbb P}^3$
and $X$ has its node
at the vertex $q$ of the cone $Q_q$, the one-dimensional family of
the ruling of $Q_q$ cuts out on $X$ a "generalized" $\g^1_3$
(i.e. an element of $\overline{W}^1_3(X)$, cf. \eqref{eq:aio1}),
with the node $q$ as a fixed point. Furthermore, the projection of
$X$ from $q$ maps $X$ $2:1$ onto the base of the quadric cone
in $H_q$. By Proposition \ref{prop:projgeom}, this determines a
$\g^1_2$ on the normalization $C$, which is a smooth,
hyperelliptic curve of genus $g=3$.  Since this holds for a general $q\in S$,
we have
$\dim (V_{|H|,1}^{hyper}) \geq 2$.

Conversely, let $X = H_X \cap Q_2 \cap Q_3$, which is irreducible,
$1$-nodal at a point $q$ and where $H_X \cong {\mathbb P}^3$.
Thus, $ T_q (S) \subset H_X$. Since $X \subset H_X \cap Q_2$, then
$X$ is contained in a quadric $Q_X$ of $H_X \cong {\mathbb P}^3$.
If we assume that $X$ has normalization $C$ which is a
hyperelliptic curve of genus $g=3$, we claim that $Q_X$ is forced
to be a quadric cone with vertex at the node $q$ of $X$. To prove
this, observe that, since $|\omega_C|$ is a special $\g^2_4$ on
$C$, then it is composed with the hyperelliptic involution on $C$.
Therefore, the map $\Phi_{|\omega_C|}$ given by the canonical
linear system is a ($2:1$)-covering of a plane conic $\Gamma$. On
the other hand, if $\nu: C \to X \subset H_X \cong  {\mathbb P}^3$
denotes the normalization morphism and if $\pi_q$ is the
projection of $X$ from its node to $ {\mathbb P}^2$, then we have
$\Phi_{|\omega_C|} = \pi_q \circ \nu$. For each $x \in \Gamma$,
let $L_{q,x} $ be the line in ${\mathbb P}^3$ connecting $x$ with
the node $q$ of $X$. Then, one can easily show that $X$ is
contained in the quadric cone ${\mathcal C}_X \subset {\mathbb
P}^3$ determined by ${\mathcal C}_X := \bigcup_{x \in \Gamma}
L_{x,q}$ having vertex at the node $q$ of $X$. Therefore,
${\mathcal C}_X = Q_X$. Thus, $\dim (V_{|H|,1}^{hyper}) = \dim (S)
= 2$ whereas $\dim(V_{|H|,1}) = 3$.

Note that  $C$  is not general from the classical Brill-Noether
theoretical point of view, indeed $\rho (3, 1, 2) = -1$, whereas
$\rho_{sing}(4,1,2,3) = 0$ and we actually have a $\g^1_2$ on $C$.

In particular, this example shows the sharpness of the bound
\eqref{eq:bound} in Corollary \ref{cor:vhdh}, whence also of Theorem 1.
}
\end{example}

\begin{remark}\label{rem:several2}\normalfont{
 Following Voisin (see \cite[pp.~ 366-367]{V1}),
Example \ref{exm1} can also be obtained as follows. Let $S$ be a
$K3$ surface with $\Pic(S) \cong \Zz[L]$, such that $L^2 = 2p-2
\geq 2$. The smooth members $D$ of $|L|$ are of geometric genus
$p$ and Brill-Noether general by Lazarsfeld's theorem quoted
in the introduction. In particular, they have the same gonality as
a generic curve.

Therefore, if $p = 2k$, $k \geq 2$, a general member $D$ is
$(k+1)$-gonal. Let $\Aa$ be the $\g^1_{k+1}$ on $D$. Following
\cite{L}, there exists a rank-two vector bundle $E:= E(D,\Aa)$ on
$S$, such that
$$
c_1(E) = L, \; c_2(E) = k+1 \; {\rm and} \; h^0(S, E) = k+2.
$$
The zero-locus of a general section of $E$ is a general member of
a $\g^1_{k+1}$ on some $D$. If $x_1, \ldots , x_{\delta}$ are
general points on $S$ such that $ 2 \delta \leq k$, the vector
space
$$
 H:= H^0 (S, \Ii_{\{x_1, \ldots, x_{\delta}\}} \otimes E)
$$
has dimension at least two. For $\alpha$ and $\beta$ general
in $H$, the curve $X:= V( \alpha \wedge \beta)$, where $\alpha
\wedge \beta \in H^0 (S, L)$, is irreducible, with nodes exactly
at the points $x_i$'s, $ 1 \leq i \leq \delta$. On the other hand,
the given two sections generate a rank-one subsheaf of the
restriction $E|_X$. If $\nu : C \to X$ is the normalization, this
rank-one subsheaf induces a line sub-bundle ${\mathcal L}_{\Aa}
\subset \nu^* (E)$ with two sections and whose moving part on $C$
is of degree $k+1-\delta$. Hence $C$ is $(k+1 - \delta)$-gonal.

If we want $C$ to be hyperelliptic, we set $\delta = k-1$. Since $2 \delta \leq k$, then
$k \leq 2$.
Therefore, the only case one obtains using this procedure, is the one discussed in Example \ref{exm1}, with
$k=2$, $p=4$ and $\delta =1$.
}
\end{remark}

The second example we want to consider is the following.

\begin{example} \label{exm2}
\normalfont{Let $p=5$ and let $S $ be a general, primitive
complete intersection surface of type $(2,2,2)$ in $\Pp^5$. We
consider the case $\delta =2$ which, by \eqref{eq:nontriv}, is the
only possible, nontrivial example to construct. By what
recalled in Remark \ref{rem:sevvar}, $V_{|H|,2} \neq \emptyset$.
Let $[X] \in V_{|H|,2}$ and let $N = \Sing(X)$;
let $\Lambda_N := l$ be the line spanned by the two nodes of $X$.
Since $C$, the normalization of $X$, has genus $g=3$ then, by
using Proposition \ref{prop:projgeom}, we first show the
following:

\vskip 7 pt

\noindent
{\bf Step 1}: {\em $\tau_l(X) = \Gamma \subset \Pp^2$ is
a conic $\Leftrightarrow$ $X \subset H = \Pp^4$ lies on a
quadric cone, which is a rank-$3$ quadric in $\Pp^4$ with vertex the line $l$ and base the conic $\Gamma$
$\Leftrightarrow$ $X$ is in the join variety of $\Gamma$ and $l$ in $\Pp^4$.}

\vskip 7pt

\noindent {\em Proof.} ($ \Rightarrow$) If $U \in |\Ii_{l/H}(1)|$,
then $\# U \cap X = 8$ and $4$ intersections are absorbed by the
two nodes. Let $\Pp^2_{\Gamma}$ be the plane containing $\Gamma$;
since $U \cap \Pp^2_{\Gamma}$ is a line, then $U \cap \Gamma$ are
two points $\gamma_1, \; \gamma_2 \in \Pp^2_{\Gamma}$. By
assumption, there exist points $n_i, \; q_i \in X$, $1 \leq i \leq
2$, such that$$\tau_X(n_i) = \tau_X(q_i) = \gamma_i, \; 1 \leq i
\leq 2;$$furthermore, the configuration of points $\{n_i, \; q_i,
\gamma_i \}$ lies on the same plane $\pi_i \subset U$ through $l$,
$1 \leq i \leq 2$. Therefore, the configuration of points $\{n_1,
\; q_1, \gamma_1, \; n_2, \; q_2, \gamma_2 \}$ lies on the
singular quadric $\Sigma_l = \pi_1 \cup \pi_2 \subset U$. If we
take the pencil of $\Pp^3$'s in $H$, passing through $l$ and
cutting 2 points on $\Gamma$, we get that $X$ is contained in the
cone over $\Gamma$ in $\Pp^4$, with vertex $l$.

\noindent
($\Leftarrow$) Let ${\mathcal Q}_l$ be a rank-$3$ quadric in $\Pp^4$ and assume that
$X \subset {\mathcal Q}_l$. Let $U \in |\Ii_{l/H}(1)|$; since
$U \cap X \subset U \cap {\mathcal Q}_l = \Sigma_l = \pi_1 \cup \pi_2$ and since
$\# U \cap X = 8$, by symmetry we have two smooth points of $X$ on each $\pi_i$,
$1 \leq i \leq 2$. By definition of the maps $\tau_X$ and $\pi_l$, we conclude.
\hfill $\Box$

Now, we want to lift-up to $S$ the geometric information proved in Step 1. Precisely, we show:

\vskip 7 pt

\noindent
{\bf Step 2}: {\em $X \subset H = \Pp^4$ lies on a rank-$3$ quadric in $H$, with the two nodes on the
linear vertex $l$ of the quadric
$\Leftrightarrow$ the surface $S = \bigcap_{i=1}^3 Q_i$ lies on a rank-$5$ quadric in $\Pp^5$
having $l$ as a generatrix, which is a $2$-secant of $S$.}

\vskip 7pt

\noindent
{\em Proof.} ($\Leftarrow$) After a possible change of basis, one can assume
that $Q_1$ is a rank-$5$ quadric in $\Pp^5$, which is a cone of vertex a point $v$
and base a smooth quadric ${\mathcal Q}$ in $\Pp^4$. If we project from $v$ to
$\Pp^4$, then $Q_1$ projects onto ${\mathcal Q}$ and $S$ onto a birationally equivalent
surface $F \subset {\mathcal Q}$, having a double curve in its singular locus. This means
that, through $v$, we have a pencil of generatrices of $Q_1$ which are
$2$-secants of $S$. Let $l$ be one of such generatrices. Let $H \in |\Ii_{l/\Pp^5}(1)|$
which is tangent to $Q_1$ along $l$. Then $H\cap S = X$ is a curve with two nodes
at the points $S \cap l$. Since $S \subset Q_1$, then $X \subset Q_1 \cap H := {\mathcal Q}_l$,
which is a rank $3$ quadric in $\Pp^4$ having $l$ as the linear vertex.

\noindent ($ \Rightarrow$) Let ${\mathcal Q}_l$ be a rank $3$
quadric in $\Pp^4$ vertex a line $l$ and assume that $X \subset
{\mathcal Q}_l$. By assumption $\Pic(S) = \ZZ[H]$, therefore $l$
is not contained in $S$. Since $X \sim H$ on $S$, the hyperplane
$H$ must be tangent to $S$ only at the two nodes of $X$, which are
$l \cap S$. Since ${\mathcal Q}_l \subset \Pp^4$ contains both $ X
= S \cap H $ and $l \subset H$ and since $S$ is non-degenerate,
irreducible and linearly normal, then ${\mathcal Q}_l$ extends to
a unique, irreducible quadric $Q_l$ in $\Pp^5$ containing $S$ and
obviously $l$ (see e.g. \cite[Lemma 7.9]{SD}, whose arguments
easily extend to the case of an irreducible, nodal hyperplane
section of $S$).

Now, $l = \Sing({\mathcal Q}_l) = \Sing(Q_l \cap H)$; therefore,
by standard facts on hyperquadrics (see e.g. \cite{Harris}, page
283), we have that either
\begin{itemize}
\item[(i)] $Q_l$ is a rank $3$ quadric in $\Pp^5$, with vertex a plane
$\Pi_l$ containing $l$, or
\item[(ii)] $Q_l$ is a rank $4$ quadric in $\Pp^5$ with vertex $l$, or
\item[(iii)] $l$ is not contained in $\Sing(Q_l)$ and $H$ is a tangent hyperplane
to $Q_l$ along $l$.
\end{itemize}

We claim that only case (iii) can occur. Indeed, in case (i), since $S \subset Q_l$ is of type $(2,2,2)$ in $\Pp^5$,
then other two linearly independent quadrics would intersect $\Pi_l$, so
$S \cap \Pi_l \neq \emptyset$, which would give that $S$ is singular, a contradiction. In case
(ii), since $S \cap l = \Sing(X)$, then $S$ would be singular too. Therefore we must be in case (iii).

Now, in case (iii), $Q_l$ must be necessarily a rank $5$ quadric with vertex a point and
$l$ a generatrix. Indeed, if $Q_l$ were smooth, then the Gaussian map would give an
isomorphism between $Q_l$ and its dual hypersurface, which is still a hyperquadric, so
the points of $l$ could not have the same tangent hyperplane.
\hfill $\Box$

\medskip

To sum up, in order to construct examples of $2$-nodal curves which are hyperplane
sections on a general, primitive, complete intersection surface $S$ of type
$(2,2,2)$ in $\Pp^5$, it is sufficient to consider the quadric cones, with vertex a point,
containing $S$ and hyperplanes which are tangent to these cones along generatrices
which are $2$-secants of $S$. Observe that such a $S$ is contained in a two-dimensional family of
hyperquadrics in $\Pp^5$. Among these hyperquadrics, we have a one-dimensional family
of quadric cones having  a zero-dimensional
vertex. Each such a cone has a one-dimensional family of $2$-secants of $S$
through its vertex. Therefore, by recalling notation as in Definition \ref{def:sevhyp},
we have, a priori, $\dim (V_{|H|,2}^{hyper} ) \geq 2.$ Furthermore,$$V_{|H|,2}^{hyper}
\subseteq V_{|H|,2}$$and, by Remark \ref{rem:sevvar},
$\dim(  V_{|H|,2} ) = 3$.

We claim that the above inclusion is proper. Indeed, let $q_1 \neq
q_2 \in S$ be such that the general element $[X] \in V_{|H|,2}$
has $\Sing(X) = \{q_1, q_2 \}$. Therefore, there exists a
hyperplane $H_X \subset \Pp^5$ s.t.$$T_{q_1}(S), \; T_{q_2}(S)
\subset H_X \;\; {\rm and} \;\; X = H_X \cap S.$$Since $H_X$
intersects each plane of $\Pp^5$ along (at least) a line, in order
to contain $T_{q_i}(S)$, $ 1 \leq i \leq 2$, we have to impose two
independent conditions on the linear system $|\Oc_{\Pp^5}(1)|$.

On the other hand, $S$ is a complete intersection of quadrics in
$\Pp^5$, i.e. $S = V(F_1, F_2, F_3)$, with $\deg(F_{i})=2$.
Therefore, to impose to a general hyperplane
$\Sigma_{i=0}^5 a_i x_i = 0$ to be tangent at $q_1$ and at $q_2$
to one of the three quadrics $V(F_i)$ gives more than two
independent conditions on the coefficients $\{a_0, \ldots, a_5\}$.

In conclusion, we have $\dim (V_{|H|,2}^{hyper} ) = 2$.

This shows that all the elements of $V_{|H|,2}^{hyper} $ on $S$
are obtained from the above construction; furthermore, this gives
examples of 2-hyperelliptic, nodal curves in the range
$\rho_{sing} > 0$. Indeed, in this case we have $g=3$, so
$\rho_{sing}(5,1,2,3) = \rho(3,1,2) + 2 = 1$. }
\end{example}

\begin{remark}\label{rem:exa}
\normalfont{In both examples presented above, $V^{hyper}_{|L|,
\delta}$ is nonempty and two-dimensional.  We will return in a
forthcoming work to the study of the nonemptyness and the
dimension of $V^{hyper}_{|L|, \delta} (S)$, for
$\rho_{sing}(p_a(L),1,2,p_a(L) - \delta)\geq 0$ and $(S,L)$ a
general, primitive $K3$.

Note also that, in Examples \ref{exm1} and \ref{exm2}, $V^{hyper}_{|L|, \delta}$ is
of codimension one in $V_{|L|, \delta}$, for $\delta= 1,\ 2$, respectively.
This shows that the property of having hyperelliptic normalizations
is not constant among all the nodal curves in a Severi variety, contrary to the case of {\em smooth} curves
(cf. the theorems of Saint-Donat \cite{SD} and Green and Lazarsfeld \cite{GL} stated in the introduction).
Moreover, the ``unexpected'' property of having hyperelliptic normalization
is not induced by a line bundle on the surface.
This indicates that the Brill-Noether theory for singular curves on $K3$ surfaces appears to be
more subtle than the one for smooth curves.
}
\end{remark}

%
\section {Neutral and non-neutral nodes} \label{S:nodchar}
%

In this section we briefly make some comments on the behavior of the nodes of a
$\delta$-hyperelliptic curve under the hyperelliptic involution, first
in the general case, and then when the curve lies
on a $K3$ surface.

Let  $X$ be a nodal curve (not necessarily lying on a $K3$
surface) and $N:= \Sing(X)$ be its scheme of nodes. If $\Aa$ is a
rank-one, torsion-free sheaf on $X$ then, by the hypotheses, one
has:
\[
\Sing(\Aa) :=
\{ q \in X \; | \; \Aa_q \; \mbox{is not locally free} \} \subseteq N \;
\mbox{and} \; \Aa_q \cong \mathfrak{m}_q \; \mbox{if} \; q \in \Sing(\Aa),
\]
where $\mathfrak{m}_q \subset \Oc_{X,q}$ is the maximal ideal (cf.
e.g. \cite{GK} or \cite{Ses}, pp. 163-165).

From now on in this section, unless otherwise stated, we shall focus on the following
situation:
\begin{notation}\label{not:node}
\normalfont{
Assume that:
\begin{itemize}
\item[(i)] $X$ is an irreducible curve of arithmetic genus $p = p_a(X)$ with $\delta$ nodes
as the only singularities;
\item[(ii)] $N = \Sing(X)$ is the (reduced) scheme of nodes of $X$;
\item[(iii)] for each subscheme $Z \subseteq N$, denote by
$\nu_Z : X_Z \to X$ the partial normalization of $X$ along $Z$. In particular,
when $Z =N$, then $X_N = C$ and $\nu_N = \nu : C \to X$ is the (total)
normalization of $X$; in this case the smooth curve $C$ has (geometric) genus
$g = p - \delta$, which is assumed to be $g \geq 2$;
\item[(iv)] there exists a line bundle $\Delta_C$ on $C$ such that
$|\Delta_C| = \g^1_2$, i.e. $C$ is hyperelliptic. From the assumption
$g \geq 2$, $\Delta_C$ is unique.
\end{itemize} }
\end{notation}

Note that the sheaf $\nu_*(\Delta_C)$ on $X$ has the following properties:
\begin{itemize}
\item $\nu_*(\Delta_C)$ is torsion-free of rank one,
\item $\deg(\nu_*(\Delta_C) ) = 2+\delta$,
\item $\Sing(\nu_*(\Delta_C)) = N$, and
\item $h^0(X, \nu_*(\Delta_C)) = h^0(C, \Delta_C) = 2$.
\end{itemize}

We now  recall the classical notion of {\em neutral} (resp.
{\em non-neutral}) pairs with respect to the $\g^1_2$ on $C$.

\begin{definition}\label{def:neutral}
Let $n \in N$ and $\nu^{-1}(n) = \{ n^{'}, n^{''} \} \subset C$. The pair $(n^{'}, n^{''})$ is called a
{\em neutral pair}
if $n^{'} + n^{''} \in |\Delta_C|$ (otherwise
the pair is called {\em non-neutral}). For brevity, the node
$n \in X$ will be called a {\em neutral node} ({\em non-neutral node}, respectively).
\end{definition}

We extend the notion of neutral (non-neutral, resp.) nodes to
linear series composed with the hyperelliptic involution of $C$.
The typical situation is given by the {\em canonical linear
system} $|\omega_C| = \g^{g-1}_{2g-2}$. Since $C$ is by assumption
 hyperelliptic, with $g \geq 2$, then
$|\omega_C| = (g-1) |\Delta_C|$.

If $|L_C| = \g^r_{2r} = \g^r_r \circ \g^1_2 =r|\Delta_C|$, then the associated morphism
\begin{equation}\label{eq:12lu12}
C \stackrel{\Psi}{\longrightarrow} \Gamma_r \subset \Pp^r,
\end{equation}
which is $2:1$ onto  a rational
normal curve $\Gamma_r \subset \Pp^r$, is the composition of the $2:1$ map given by $|\Delta_C|$, and of the
$r$-tuple Veronese embedding. Hence,  it makes sense to give the following:

\begin{definition}\label{def:12lu13}
Let $n \in N$ be a node of $X$ and let $(n^{'}, n^{''})$ be the
corresponding pair on $C$. Then $(n^{'}, n^{''})$ is a {\em
neutral pair} (equiv., $n$ is a {\em neutral node}) for $|L_C|$ if
and only if it is a neutral pair (equiv., a neutral node) for
$|\Delta_C|$.
\end{definition}

\begin{proposition}\label{prop:12lu14}
Let $n \in N$ and let $|L_C| = \g^r_{2r}$, with $r \leq g-1$, be a
base point free linear system on $C$ which is composed with
$|\Delta_C|$. Then $\nu_*(L_C)$, as a rank-one torsion-free sheaf
on $X$, is not globally generated at the neutral nodes, and
globally generated elsewhere.
\end{proposition}
\begin{proof} First of all, observe that as a straightforward
consequence of Definition \ref{prop:12lu14}, $\nu_*(L_C)$ is not
g.g. at a point $p \in X$ if and only if $\nu_*(\Delta_C)$ is.
Therefore,  since $\nu$ is an isomorphism outside $N$, it is
sufficient to study the generation of $\nu_*(\Delta_C)$ at the
nodes of $X$. Using Notation \ref{not:node}, let $A\subseteq N$ be
the (possibly empty) subset of neutral nodes and $Z\subseteq N$ be
the (possibly empty) subset of non-neutral nodes of $X$. Consider
the following commutative diagram:
\begin{equation}\label{eq:12lu2}
 \xymatrix{
     C  \ar[dr]_{\nu} \ar[r]^{\tau} & X_Z \ar[d]^{\varphi} & \hspace{-1cm}   \\
     & X, &
 }
 \end{equation}
where
$\varphi$ normalizes the nodes in $Z$ (so $X_Z$ has still
$a:= \Card (A)$ nodes),
and
$\tau=\nu_A$ normalizes the residual $a$ nodes. Denote by
$f:= f_{\Delta_C} : C \to \Pp^1$
the ($2:1$)-morphism induced by $|\Delta_C|$.

Since the map $\tau$ identifies $a$ neutral pairs, we have that the morphism $f$ induces on $X_Z$ a ($2:1$)-morphism
$h: X_Z \to \Pp^1$,
which gives the commutative diagram
\begin{equation}\label{eq:12lu6}
 \xymatrix{
     C  \ar[dr]_{f} \ar[r]^{\tau} & X_Z \ar[d]^{h}   \\
     &  \Pp^1,
 }
 \end{equation}
i.e. $\Delta_C$ induces a line bundle $R \in \Pic^2(X_Z)$ such that:
\begin{itemize}
\item[(i)] $R$ is globally generated on $X_Z$ (since $h$ is a morphism), and
\item[(ii)] $\tau^*(R) = \Delta_C$ (since the $\g^1_{2}$ on $C$ is unique).
\end{itemize}

From (ii), we have that $\tau_{*}(\Delta_{C})=R\otimes \tau_{*}\Oc_{C}$. From the exact sequence
$$
 0\to \Oc_{X_Z}\to \tau_{*}\Oc_{C}\to \Oc_{A} \to 0
$$ tensored by $R$, we see that the cokernel of the evaluation
morphism$$ H^0 (\tau_{*}(\Delta_C)) \otimes \Oc_{X_Z}\cong
H^0(R)\otimes \Oc_{X_Z} \to \tau_{*}(\Delta_C) $$ is exactly
$\Oc_{A}$.  So  $\nu_{*}(\Delta_{C})$ is not generated at the
neutral nodes.

On the other hand, if  $\varphi_*(R)$ were not globally generated
as a rank-one torsion-free sheaf on $X$, then consider $\La :=
\im(ev)$, where $ev: H^0(\varphi_*(R))\otimes \Oc_{X} \to
\varphi_*(R)$ is the (non-surjective) evaluation map. It is easy
to show that $\La\in \Pic (X)$ and $\varphi^*(\La)= R$. The latter
is impossible, since $\varphi: X_Z \to X$ identifies the non-neutral
pairs, so $\varphi_*(R)$ is globally generated.

From Diagram \eqref{eq:12lu6}, it follows that $\nu_*(\Delta_C)$,
as rank-one torsion-free sheaf on $X$, is globally generated at
the non-neutral nodes of $X$. This concludes the proof.
\end{proof}

\vskip 10pt

\begin{remark} \label{rmk:nnchar}
{\rm We note that from the last proposition and the discussion above, it is clear
that the non-neutral nodes $N_0$ of $X$ can be characterized by any of the three
following equivalent descriptions}:
\begin{itemize}
\item[(i)]   $N_0$ is the set of the nodes whose two preimages on the normalization of $X$ are
not conjugated by the hyperelliptic involution.
\item[(ii)] $N_0$ is the minimal set of nodes such that the curve obtained by desingularizing these
nodes possesses a $\g^1_2$, i.e. it admits a $(2:1)$-morphism to $\Pp^1$.
\item[(iii)] $N_0$ is the set of nodes at which $\nu_*(\omega_C)$ is globally generated.
\end{itemize}
{\rm If we furthermore assume that $X$ lies on a $K3$ surface and
$X \in |L|$ for $L \in \Pic (S)$,
then by \eqref{eq:9lu*}
and by the short exact sequence
\[ 0 \to \Oc_S \to \Ii_{N/S} \* L \to \Ii_{N/X} \* L \to 0, \]
the conditions (i)-(iii) are equivalent to:}
\begin{itemize}
\item[(iv)] $N_0$ is the set of nodes at which $\Ii_{N/S} \* L$ is
globally generated as a sheaf on $S$.
\end{itemize}
\end{remark}

\vskip 10pt

If $|L|$ is hyperelliptic on a $K3$ surface $S$, then $\varphi_L$ maps the surface $S$
generically $2:1$ onto a rational ruled surface or onto the Veronese surface in $\Pp^5$, by
Saint-Donat \cite{SD}, so that any irreducible member of $|L|$ is mapped $2:1$
onto a smooth rational curve. Hence, if $V_{|L|, \delta}$ is not
empty, any node of a curve in it must be neutral. To see such an example,
it is sufficient to consider $S$ as the general degree-$2$ $K3$ and take $X$
as the pull-back of a line tangent to the branched plane sextic curve.

\medskip

The next two examples show that the nodes in Examples \ref{exm1} and
\ref{exm2} above are all non-neutral. This will also follow from the
results in the next section (more precisely from Theorem 2). However, since this can be
seen geometrically, we find it instructive to include these examples at this point.

\begin{example}\label{exa:1} {\normalfont
Let $S \subset \Pp^4$ be a general complete intersection of type
$(2,3)$ and let $H$ be the hyperplane section of $S$. From Example
\ref{exm1}, we know that $V_{|H|,1}^{hyper} \neq \emptyset$. Let
$[X] \in V_{|H|,1}^{hyper}$ and let $\{p\} = Sing (X)$. Let $\nu :
C \to X$ be the normalization morphism and let $\nu^{-1}(p) =
\{p', p''\}$. Suppose that $p$ is a neutral node of $X$, i.e.
$|\Oc_C(p' + p'')| = |\Delta_C | = \g^1_2$ on $C$. Thus,
\begin{equation}\label{eq:aq+}
h^0(C, \omega_C \otimes \Oc_C(- p' - p'')) = h^0(C, \omega_C
\otimes \Oc_C(- p'))=2.
\end{equation}On the other hand, by \eqref{eq:9lu1}, formula
\eqref{eq:aq+} becomes
\begin{equation}\label{eq:aq++}
h^0(C, \nu^*(H)- 2p' - 2 p'') = h^0(C, \nu^*(H)- 2p' - p'')=2,
\end{equation}which is a contradiction.
Indeed, the geometric meaning of \eqref{eq:aq++} is that there
should exist a pencil of planes in $\Pp^3=H$ such that each plane
of the pencil passes through $p$ and, since it is tangent at $p$
to one of the two branches of $X$ around $p$, it must be tangent
also to the other branch of $X$ through $p$. This cannot occur
since only one such plane does exist; thus $p$ is a non-neutral
node for $X$. }
\end{example}

\begin{example}\label{exa:2}{\normalfont From Example \ref{exm2}, let
$S \subset \Pp^5$ be a general, primitive $K3$ surface of type
$(2,2,2)$ and let $[X] \in V_{|H|,2}^{hyper}(S)$. Let $\Sing(X) =
\{ n,q \}$. We want to show that both $n$ and $q$ are non-neutral
nodes for $X$.

By contradiction, two cases must be considered. Let $\nu : C \to
X$ be the normalization. Then:

\noindent
{\bf Case 1:} Suppose $n$ neutral and $q$ non-neutral.
Let $H_X \subset \Pp^5$ be the hyperplane such that $X= S \cap
H_X$.

If $\nu^{-1}(n) = \{n', n'' \}$ and $\nu^{-1}(q) = \{q', q'' \}$,
then $|\Oc_C(n' + n'')| = |\Delta_C | = \g^1_2$ on $C$. Thus, as
in Example \ref{exa:1}, we would have:
\begin{equation}\label{eq:aq**}
h^0(\nu^*(H)- 2n' - 2 n'' - q' - q'') = h^0(\nu^*(H)- 2n' - n'' -
q' - q'')=2,
\end{equation}which is a contradiction. Indeed, as above, \eqref{eq:aq**} would imply
there exist a pencil of hyperplanes in $\Pp^4 = H_X$ such that
each hyperplane of the pencil passes through $q$ and, since it is
tangent at $n$ to one of the two branches of $X$ through $n$, it
must be tangent also to the other branch of $X$ through $n$; in
fact only one such hyperplane actually exists in $\Pp^4 = H_X$.

\noindent
{\bf Case 2:} Suppose both $n$ and $q$ neutral. One can
conclude as above. }
\end{example}

%
\section {$\delta$-hyperelliptic nodal curves: a
classification result for $\rho_{sing}<0$.}\label{S:BR}
%

As already mentioned in the introduction, here we want to classify
linear systems $|L|$ on a $K3$ surface
$S$ with {\em arbitrary} Picard group, such that, as in Definition \ref{def:sevhyp},
$V_{|L|, \delta}^{hyper} \neq \emptyset$
and $\rho_{sing}<0$.
This can be viewed as a natural extension to the nodal case of
Saint-Donat's result.

In the particular case of $S$ with $\Pic(S) = \ZZ[L]$ and $p =
p_a(L) \geq 3$, we find again Corollary \ref{cor:vhdh} via a
different approach.

The techniques used here are related to a Bogomolov-Reider type approach
for separation of points by the linear system $|L|$ on $S$.

By recalling Notation \ref{not:node}, we first prove the following
more general result, which translates the property of having
k-gonal normalization into a failure of separation of
zero-dimensional schemes by $|L|$.

\begin{theorem}\label{thm:hyper}
Let $S$ be a $K3$ surface and $L$ be a globally generated line
bundle on $S$ of sectional genus $p =p_a(L) \geq 3$. Let $\delta$
be an integer such that $\delta \leq p-2$. Let $[X] \in V_{|L|,
\delta}$ and let $N_0  \subseteq \Sing (X):=N$ be a subset of the
scheme of the $\delta$ nodes of $X$. Let $X_{N_0}$ be the curve
obtained by desingularizing $X$ at the nodes in $N_0$.

Then $X_{N_0}$ possesses a base point free, complete $\g^1_k$, for
$k \geq 2$, if and only if there exist $k$ points $q_1, \ldots,
q_k \in X \setminus N$ such that the zero-dimensional subscheme
$N' : = N_0 \cup \{ q_1, \ldots, q_k \} \sub X \sub S$ does not
impose independent conditions on the linear system $|L|$, but for
any proper subset $Z \subsetneq \{q_1, \ldots, q_k\}$, the
subscheme $N_0 \cup Z$ imposes independent conditions on $|L|$;
more precisely
\begin{equation} \label{eq:P}
  h^1(\Ii_{N'/S} \otimes L ) = 1 \; \mbox{and} \; h^1(\Ii_{(N_0 \cup Z) /S} \otimes L )
  = 0 \; \mbox{for any} \; Z \subsetneq \{q_1, \ldots, q_k\}
\end{equation}

Moreover, this property holds for every $k$-uple of points $q_1,
\ldots, q_k \in X \setminus N$ (not necessarly distinct) such that
their inverse images on $X_{N_0}$ are part of the
$\mathfrak{g}^1_k$.
\end{theorem}

\begin{proof}
Denote by $\mu_{N_0} : S_{N_0} \to S$ the blow-up of $S$ along
${N_0}$, which induces the morphism $\nu_{N_0} : X_{N_0} \to X$,
which normalizes the nodes in $N_0$. Let $B := \sum_{i=1}^{\delta}
E_i$ be the total $\mu_{N_0}$-exceptional divisor and let
$\tilde{N}_0 = \Supp(\Oc_{X_{N_0}}(B)) $ be the scheme of
$2\delta_0$ points of $X_{N_0}$, the pre-images of the $\delta_0
\leq \delta$ nodes in $N_0$. As in \eqref{eq:9lu1},
$\omega_{X_{N_0}} \cong \Oc_{X_{N_0}}( \nu^* L) \otimes
\Oc_{X_{N_0}} (- \tilde{{N_0}})$, where $\omega_{X_{N_0}}$ is the
dualizing sheaf of ${X_{N_0}}$. This gives
\begin{equation}\label{eq:tangsev}
H^0 ({X_{N_0}} , \omega_{X_{N_0}}) \cong \frac{H^0 (S,
\Ii_{{N_0}/S} (L))}{H^0(\Oc_S)}.
\end{equation}

Let $q_1', \ldots, q_k' \in {X_{N_0}} - \tilde{{N_0}} -
\Sing({X_{N_0}})$ be $k$ (not necessarily distinct) points and
$q_1, \ldots, q_k \in X \setminus N$ their images via $\nu$.
Consider $\pi$ the further blow-up of $S_{N_0}$ along $\{q_1',
\ldots, q_k'\}$. Therefore, $\mu_{N'} :=  \pi \circ \mu_{N_0} $ is
the total blow-up of $S$ along $N' := {N_0} \cup \{q_1, \ldots,
q_k \}$. Denote by $E_{q_1}, \ldots, E_{q_k}$ the
$\mu_{N'}$-exceptional divisors over $q_1, \ldots ,q_k$,
respectively.

Denote by $\overline{X}_{N_0}$ the proper transform of $X$ in
$S_{N'}$. Then ${X_{N_0}} \cong \overline{X}_{N_0}$; let $q_j''$
be the point on $\overline{X}_{N_0}$ corresponding to $q_j'$ on
$X_{N_0}$, $ 1 \leq j \leq k$. Consider the standard exact
sequences
\begin{equation}\label{eq:blowup1}
0 \to  \Oc_{S_{N'}}(B)  \to  \Oc_{S_{N'}}(\mu_{N'}^*(L) - B -
E_{q_1} - \cdots - E_{q_k})
 \to \omega_{\overline{X}_{N_0}} (-q_1'' - \cdots - q_k'')  \to 0
\end{equation}
and
\begin{equation}\label{eq:blowup2}
0 \to  \Oc_{S_{N'}}(B + E_{q_1} + \cdots + E_{q_k})  \to
\Oc_{S_{N'}}(\mu_{N'}^*(L) - B )
 \to  \omega_{\overline{X}_{N_0}}   \to 0
\end{equation}
on $S_{N'}$. By the Leray isomorphism, Fujita's Lemma (cf.
\cite{KMM}) and by the fact that $S$ is $K3$, we get that $h^1
(\Oc_{S_{N'}}(B)) = h^2 (\Oc_{S_{N'}}(\mu_{N'}^*(L) - B - E_{q_1}
- \cdots - E_{q_k}))= 0$ and $h^2(\Oc_{S_{N'}}(B))=1$. Hence, by
the Leray isomorphism again, we have
\begin{equation} \label{eq:tangsev2}
H^0 (\omega_{\overline{X}_{N_0}}) \cong \frac{H^0 (S,
\Ii_{{N_0}/S} (L))}{H^0(\Oc_S)}, \;\; H^0
(\omega_{\overline{X}_{N_0}} (-q_1'' - \cdots - q_k'') ) \cong
\frac{H^0 (S, \Ii_{N'/S} (L))}{H^0(\Oc_S)},
\end{equation}
and
\begin{equation} \label{eq:tangsev2'}
h^1 (\omega_{\overline{X}_{N_0}} (-q_1'' - \cdots - q_k'') ) = h^1
(S, \Ii_{N'/S} (L)) - 1.
\end{equation}

From Serre duality and \eqref{eq:tangsev2'} we see that
$|\Oc_{\overline{X}_{N_0}} (q_1' + \cdots + q_k')|$ is a $\g^1_k$
if and only if $h^1(\Ii_{N'/S} (L))=1$ or, equivalently, from
\eqref{eq:tangsev2} if and only if $N'$ imposes only
$\delta_0+k-1$ conditions on $|L|$.

At the same time, the $\g^1_k$ is base point free if and only if,
for any proper subset $Z \subsetneq \{q_1, \ldots, q_k\}$,
denoting the corresponding subscheme on $\overline{X}_{N_0}$ by
$Z''$, we have $h^0(\Oc_{\overline{X}_{N_0}} (Z''))=1$, which is
equivalent to saying that $h^1(\Ii_{(N_0 \cup Z)/S} (L))=0$, or
that $N_0 \cup Z$ imposes independent conditions on $|L|$, as
$h^1(L) = 0$ since $L$ is big and nef.

To conclude, we have left to show that if $N'$ does not impose
idependent conditions on $|L|$, but $N_0 \cup Z$ does, for any
proper subset $Z \subsetneq \{q_1, \ldots, q_k\}$, then
$h^1(\Ii_{N'/S} (L))=1$.

Assume, to get a contradiction, that $h^1(\Ii_{N'/S} (L))>1$.
Then, as in \eqref{eq:tangsev2'}, we would have that
$h^0(\Oc_{\overline{X}_{N_0}} (q_1'' + \cdots + q_k'')) >2$.
Therefore, we can find a proper subset $Z \subsetneq \{q_1,
\ldots, q_k\}$, such that, denoting the corresponding subscheme on
$\overline{X}_{N_0}$ by $Z''$, we have $h^0(\overline{X}_{N_0}
(Z'')) \geq 2$. Hence, arguing as above, $N_0 \cup Z$ does not
impose independent conditions on $|L|$, a contradiction.
\end{proof}

Using Theorem \ref{thm:hyper} and the Reider-like results on
higher order embeddings of $K3$ surfaces in \cite{K}, we are now
ready to prove Theorem 2, stated in the introduction, which
classifies linear systems $|L|$ on $S$ for which $V_{|L|,
\delta}^{hyper} \neq \emptyset$, with $\rho_{sing}<0$. We remark
that, as already mentioned in the introduction, a similar
statement can also be obtained for curves with normalizations
possessing $\g^1_k$'s for $k >2$, using Theorem \ref{thm:hyper}
and the corresponding notion of neutral and non-neutral nodes for
$\g^1_k$'s. The principle is the same, but the classification
becomes longer.

\begin{proof}[Proof of Theorem 2.] Assume that $V_{|L|,
\delta}^{hyper} \neq \emptyset$ and let $[X]$ be a point of this
scheme. Let $N_0 \subseteq N:=\Sing (X)$ be the set of non-neutral
nodes of $X$, which consists of $\delta_0 \leq \delta$ points, by
assumption. Then, by Theorem \ref{thm:hyper} and Remark
\ref{rmk:nnchar}, there exist two points $p , q \in X \setminus N$
such that, if $N' = {N_0} \cup \{p,q \}$, then $|L|$ does not
separate $N'$.

Let $Z \subseteq N'$ be minimal, so that for each $Z' \subsetneq
Z$,  $|L|$ separates $Z'$ but it does not separate $Z$. We claim
that $Z = N'$.

Indeed, if  $p \not \in Z$ then in particular $|L|$ does not
separate ${N_0} \cup \{q \}$ and by (\ref{eq:tangsev}) it follows that
the inverse image of $q$ is a base point of $|\omega_{X_{N_0}}|$, a
contradiction. By symmetry we therefore have that both $p, q \in
Z$. If $Z = N_1 \cup \{p,q \}$ for some $N_1 \subsetneq N_0$, then
already a partial normalization of $\delta_1:=\length N_1$ of the nodes of $X$
would admit a $\g^1_2$,
contradicting property (ii) in Remark \ref{rmk:nnchar}. Hence $Z = N'$.

We have that $\length (N') = \delta_0 + 2$. Assume now that
\[
\delta_0 \leq \frac{p_a(L) -3}{2},
\]
or equivalently
\begin{equation}\label{eq:aiuto1}
L^2 \geq 4 (\length (N') -1).
\end{equation}

By \cite[Theorem 1.1]{K}, the condition \eqref{eq:aiuto1} implies
that there exists an effective divisor $D$ on $S$ passing through
$N'$ and satisfying the numerical conditions
\begin{eqnarray}
\nonumber & 2D^2 \stackrel{(i)}{\leq} L.D \leq D^2 + \delta_0+2
\stackrel{(ii)}{\leq}
2\delta_0+4   \\
\label{eq:kva} & \mbox{ with equality in (i) if and only if } L
\sim 2D \mbox{ and }
L^2 \leq 4\delta_0+8,   \\
& \mbox{and \nonumber equality in (ii) if and only if } L \sim 2D
\mbox{ and } L^2= 4\delta_0+8.
\end{eqnarray}

Setting $D': = L- D$ one easily finds that \eqref{eq:kva} implies
\begin{equation}
  \label{eq:kva2}
D'. D \leq \delta_0 + 2.
\end{equation}
Furthermore, by \cite[Lemma 3.6(v)]{K} (or the proof of
\cite[Prop. 1.13]{jk}), we have $L.(L-2D) \geq 0$, whence
\begin{equation}
  \label{eq:kva2'}
{D'}^2 \geq D^2.
\end{equation}
By \eqref{eq:aiuto1} and \eqref{eq:kva} we have
\[ L^2 \geq 4\delta_0 +4 > 2\delta_0 +4 \geq D.L, \]
whence $D \not \supseteq X$, by the nefness of $L$. Therefore,
since $D \supset N'$ we have by Bezout's Theorem
\begin{equation}\label{eq:kva3}
 D. X = D. L \geq 2 \delta_0 + 2.
\end{equation}
This gives
\begin{equation}
  \label{eq:kva4}
D^2 = D. L - D. D' \geq 2 \delta_0 + 2 - \delta_0 - 2 = \delta_0
>0.
\end{equation}

The Hodge index theorem and (\ref{eq:kva2}) yields
\begin{equation}\label{eq:kva5} (D^2)^2 \leq D^2 {D'}^2 \leq (D. D')^2 \leq
(\delta_0 + 2)^2,
\end{equation}
which gives $\delta_0 \leq D^2 \leq \delta_0 +2$ (note that, $D^2$
is even, since the intersection form on a $K3$ surface is even).

If $D^2=\delta_0+2$ we get from (\ref{eq:kva}) that $L \sim 2D$,
so that $L^2=4\delta_0+8$, i.e. $p_a(L)=2\delta_0+5$, which is the
second case in (ii) of the theorem.

If $D^2=\delta_0+1$ we get from (\ref{eq:kva}) that either
$L.D=2\delta_0+2$ with $L \sim 2D$, or $L.D=2\delta_0+3$.

In the first case we get $L^2=4\delta_0+4$, i.e.
$p_a(L)=2\delta_0+3$, which is the first case in (ii). Moreover,
since $L$ is nonhyperelliptic by assumption, we must have $D^2
\geq 4$, whence $\delta_0 \geq 3$, by Saint-Donat's theorem quoted
in the introduction.

In the second we get $D.D'=D.L-D^2=\delta_0+2$. From
(\ref{eq:kva2'}) and (\ref{eq:kva5}) and the fact that $\delta_0$
is odd we get that ${D'}^2=\delta_0+1$ or $\delta_0+3$, yielding
$L^2=D^2+{D'}^2+2D.D'=4\delta_0+6$ and $4\delta_0+8$ respectively.
In the first subcase we get from the results in \cite{K} (see the
precise statement in \cite[Prop. 1.13]{jk}) that $L \sim
2D+\Gamma$ for a smooth rational curve $\Gamma$ such that
$\Gamma.D=1$ and $\Gamma \cap N' \neq \emptyset$. But since
$\Gamma.L=0$ we have $\Gamma \cap X = \emptyset$, a contradiction.
In the second subcase we compute $(L-2D)^2=0$ and $(L-2D).L=2$ so
that $|L|$ is a hyperelliptic system by Saint-Donat's theorem
quoted in the introduction, a contradiction.

Finally, if $D^2=\delta_0$ then $\delta_0$ is even and we get from
(\ref{eq:kva}) and (\ref{eq:kva3}) that $L.D=2\delta_0+2$. Then
(\ref{eq:kva5}) yields ${D'}^2= \delta_0$, $\delta_0+2$ or
$\delta_0+4$ or $(\delta_0,{D'}^2)=(2,8)$. The last case is case
(i) in the theorem and one easily sees that the first three cases
yield the three cases in (iii).

We have therefore proved that if $\delta_0 \leq \frac{p_a(L)
-3}{2}$ we must be in one of the three cases (i)-(iii) of Theorem
2. In particular we have proved that either $\delta_0 \geq (p_a(L)
-5) / 2$ or $(\delta_0, p_a(L))=(2,10)$, as stated in the theorem.
\end{proof}

Now Corollary 3 in the introduction immediately follows.

\begin{proof}[Proof of Corollary 4]
If $\delta \leq \frac{p_a(nL)-3}{2}=\frac{1}{4}n^2L^2-1$, then $\delta_0 \leq  \frac{p_a(nL)-3}{2}$,
where $\delta_0$ is the number of non-neutral nodes of $X \in V_{|nL|, \delta}^{hyper}$. By
Theorem 2 and the fact that $\Pic (S) \cong \Zz[L]$, we must be in one of the cases (i) or (ii) of
Theorem 2.

In case (i) we must have $n=3$, $L^2=2$ and $2=\delta_0 \leq \delta \leq 3$, whence $\delta=2$ or
$3$. This yields case (a) in the corollary.

In case (ii), we have $nL \sim 2 D$, with $D \sim kL$ for some $k
\geq 1$, whence $n=2k$. Moreover $D^2= \delta_0+1$ or
$\delta_0+2$.

If $D^2= \delta_0+1$ we find $\delta_0 = \frac{1}{4}n^2L^2-1$, so
that $\delta=\delta_0$. It follows that $k=\sqrt{\frac{\delta +
1}{L^2}}$ and $n=2\sqrt{\frac{\delta + 1}{L^2}}$, which is case
(b) in the corollary.

If $D^2= \delta_0+2$ we find $\delta_0 = \frac{1}{4}n^2L^2-2$, so
that $\delta_0=\delta-1$ or $\delta$. In the first case we end up
in (b) as above, and in the second we get $k=\sqrt{\frac{\delta +
2}{L^2}}$ and $n=2\sqrt{\frac{\delta + 2}{L^2}}$, which is case
(c) in the corollary.
\end{proof}

We remark that for what concerns the case $\delta=1$
we have a complete picture.

\begin{corollary} \label{cor:unnodo}
Let $S$ be a $K3$ surface and $L$ be a globally generated, nonhyperelliptic
line bundle on $S$ of sectional genus $p=p_a(L) \geq 3$. Then:
\begin{itemize}
\item[(i)] $V_{|L|, 1}^{hyper} = \emptyset, \;\; \mbox{if} \;\;  p \geq 5$;
\item[(ii)] if, furthermore, $S$ is a general, primitive $K3$ such that $\Pic(S) = \Zz[L]$, then
$$\dim (V_{|L|, 1}^{hyper}) =2, \;\; \mbox{for} \;\; p=3,4.$$
\end{itemize}
\end{corollary}

\begin{proof} Statement (i) is an immediate consequence of Theorem 2.
Statement (ii) follows from Example \ref{exm1} when $p=4$. For $p=3$, just take $S$ as
general in $\Pp^3$ embedded by $|L|$ and, through the general
point $p \in S$, take a hyperplane section which is tangent to $S$
at $p$. This yields an irreducible curve on $S$ with one node,
that is a curve of geometric genus two, which is automatically
hyperelliptic. Therefore $\dim (V_{|L|, 1}^{hyper})=
\dim(V_{|L|,1}) =2$ in this case as well.
\end{proof}

%
%

\section*{Acknowledgments} We are extremely grateful to C.
Ciliberto for many valuable conversations and helpful comments on
the subject and for having pointed out a mistake in a preliminary
version of the paper. We also thank A. Verra for pointing out key
examples. We finally express our gratitude to the Department of
Mathematics of the Universit\`a "Roma Tre", where the main part of
this work has been done, for the nice and warm atmosphere as well
as for the kind hospitality.

%
%

\end{document}